\documentclass[12pt,a4paper,reqno]{amsart}
\usepackage{amsthm}
\usepackage{amssymb}
\usepackage{mathtools}
\usepackage{pdfpages}
\usepackage{changes}
\usepackage{import}
\usepackage{color} 
\definecolor{darkblue}{rgb}{0,0,0.6}
\usepackage[breaklinks, pdftex, ocgcolorlinks,colorlinks=true, citecolor=darkblue, filecolor=darkblue, linkcolor=darkblue, urlcolor=violet,  linktocpage=true, backref=page]{hyperref}

\renewcommand*{\backref}[1]{}
\renewcommand*{\backrefalt}[4]
{%
    \ifcase #1 (Not cited.)%
        \or        (Cited on page~#2.)
        \else      (Cited on pages~#2.)
    \fi
}

\newlength{\bibitemsep}
\newlength{\bibparskip}
\let\oldthebibliography\thebibliography
\renewcommand\thebibliography[1]{%
  \oldthebibliography{#1}%
  \setlength{\parskip}{\bibitemsep}%
  \setlength{\itemsep}{\bibparskip}%
}

\usepackage[colorinlistoftodos]{todonotes}

\makeatletter

\theoremstyle{plain}
\newtheorem{theorem}{Theorem}[section]
\newtheorem{proposition}[theorem]{Proposition}

\newtheorem{corollary}[theorem]{Corollary}

\theoremstyle{definition}
\newtheorem{definition}[theorem]{Definition}

\usepackage{etoolbox}
\patchcmd{\thebibliography}{\leftmargin\labelwidth }{\leftmargin\labelwidth\labelsep=20pt}{}{}

\begin{document}

\title[]{An Associative Left Brace is a Ring}

\begin{abstract}
    It is proved that if 
    a left brace $A$ has the operation $\ast$ associative, then $A$ is a two-sided brace. Consequently, $A$ is a Jacobson radical ring. 
    This answers a question of {Ced{\'o}}, {Gateva-Ivanova} and Smoktunowicz.
\end{abstract}


\author{Ivan Lau}

\address{Department of Mathematics, Simon Fraser University, 8888 University
Drive, Burnaby, BC, Canada, V5A 1S6}
\email{iplau@sfu.ca}

\maketitle

\section{Introduction}
Rump in \cite{Rump07} generalized Jacobson radical rings to a new algebraic structure, now known as left braces (see \cite{Ced18} for a survey).
The main motivation to study left braces is their relation with a particular type of solutions of the Yang-Baxter equation, the non-degenerate involutive set-theoretic solutions, in the sense of \cite[\S 3]{CJO14}.
 Following this approach, understanding the structure of non-degenerate involutive set-theoretic solutions of the YBE is reduced to classification of left braces in \cite{BCJ16}.

A left brace is a triple $(A, +, \circ)$ where $(A, +)$ is an abelian group, $(A, \circ)$ is a group and $a \circ (b + c) + a = a \circ b + a \circ c$ for all $a, b, c \in A$. If we define an operation $\ast$ on $A$ by $a \ast b = a \circ b - a - b$ for all $a, b \in A$, then the last axiom of a left brace is equivalent to $\ast$ being left-distributive with respect to the operation $+$. A Jacobson radical ring is a left brace $(R, +, \circ)$ satisfying more requirements: $(R, +, \ast)$ is a non-unital ring.
Equivalently, we say that a left brace $(A, +, \circ)$ is an algebraic structure similar to a Jacobson radical ring except $(A, +, \ast)$ need not satisfy two of the ring axioms:  right-distributivity and associativity of $\ast$. 

It may seem that the two axioms relaxed are both essential for a left brace to be a Jacobson radical ring. Rump shows that  for a left brace $A$, the operation $\ast$ being right-distributive 
is sufficent for $(A, +, \ast)$  to be a Jacobson radical ring. In other words, a left brace $A$ with the operation $\ast$ being right-distributive also has the operation $\ast$ being associative.
The proof is rather illuminating: the associativity of the operation $\ast$ is ``supplied" by the associativity of the two groups $(A, +)$ and $(A, \circ)$.

It is asked in \cite[Question 2.1(2)]{CGS18} whether the operation $\ast$  being associative is also a sufficiency condition.
We answer 
this question positively through Corollary \ref{maincor}, which follows immediately from our main result, Theorem \ref{main}. 
Prior to this, Theorem \ref{main} has been shown to hold for left
braces of odd order in \cite[Proposition 2.2]{CGS18}. 
It is also shown to hold for 8683 left braces using GAP in  \cite[\S 2.2]{KSV18}. 
As mentioned, the following is our main result:

\begin{theorem}\label{main}
    Let $(A, +, \circ)$ be  a  left  brace. Suppose that the operation $\ast$, 
    defined by $a \ast b = a \circ b - a - b$ for all $a, b \in A$ is associative. Then $A$ is  a  two-sided  brace.
\end{theorem}

Rump shows that a left brace $(A, +, \circ)$ is a two-sided brace if and only if $(A, +, \ast)$ is a Jacobson radical ring in \cite[p. 159]{Rump07} (see also \cite[Proposition 1]{CJO14}). This gives us the following corollary.

\begin{corollary}
\label{maincor}
    Let $(A, +, \circ)$ be  a  left  brace. Suppose that the operation $\ast$, 
    defined by $a \ast b = a \circ b - a - b$ for all $a, b \in A$ is associative. Then $(A, +, \ast)$ is  a Jacobson radical ring.
\end{corollary}

The proof of Theorem \ref{main} is given in Section \ref{Sect: 3}, after we review the background on left braces in Section \ref{Sect: 2}.

\section{Background on Left Braces}  \label{Sect: 2}
In this section, we will review the background on left braces, 
sufficient to discuss Theorem \ref{main}. 
Among the propositions discussed, only Proposition \ref{prop: basic circle arithmetic} is used directly in the proof of Theorem \ref{main}.

\begin{definition}
    A \textit{left brace} is a triple $(A, +, \circ)$, 
    where $(A, +)$ is an abelian group and 
    $(A, \circ)$ is a group such that 
    \[
          a \circ (b + c) + a  
        = a \circ b  + a \circ c 
        \tag{2.1} \label{eqn: left distributive}
    \]
    for all $a,b,c \in A$. 
\end{definition}

\begin{definition}
    A \textit{right brace} is defined analogously, 
    replacing condition (\ref{eqn: left distributive}) by
    \[
          (a + b) \circ c + c 
        = a \circ c  + b \circ c 
        \tag{2.2} \label{eqn: right distributive}
    \]
    for all $a,b,c \in A$. 
\end{definition}

\begin{definition}
    If $(A, +, \circ)$ is both a left and 
    a right brace (for the same operations), 
    then $A$ is called a \textit{two-sided brace}.
    \end{definition}

\begin{definition}
    For any left brace $A$, we define the operation $\ast$ by 
    \[
        a \ast b =  a \circ b - a - b 
        \tag{2.3}\label{eqn: star}
    \]
    for all $a, b \in A$. 
\end{definition}

It is straightforward to verify that the operation $\ast$ 
is left distributive with respect to the operation $+$ 
(hence the name left brace), that is,
\[
      a \ast (b + c) 
    = a \ast b + a \ast c 
    \tag{2.4}\label{eqn: left distributive for star}
\]
for all $a, b, c \in A$.
On the contrary, the operation $\ast$ is not right distributive 
nor associative for a left brace in general. 

Throughout the paper, we write 0 as the additive identity of a left brace $A$ and $-a$ as the additive inverse of $a$ for any $a \in A$. 
Let $a$ be any element of a left brace $A$ and $n \geq 0$ be an integer. 
We write $na$ as $\underbrace{a + a + \cdots + a}_{n \text{ times} }$. 
Similarly, if $m < 0$ is a negative integer, we write $ma$ to represent $\underbrace{-a -a - \cdots - a}_{-m \text{ times} }$. 
It is straightforward to check that the additive identity of $A$, 0, is also the identity of $(A, \circ)$. For any $a \in A$, we write $a^{-1}$ as the inverse of $a$ under the operation $\circ$.

We now develop some identities for left braces. 
While the proof techniques employed are standard for readers with a background in ring theory,
we will sketch the proofs for the convenience of the readers.

\begin{proposition}
\label{prop: basic star arithmetic}
    Let $(A, +, \circ)$ be a left brace. Then for all $a, b, c \in A$, we have
    \begin{itemize}
        \item [(i)] $a \ast 0 =  0 \ast a  = 0$.
        \item [(ii)] $a \ast (-b) = -(a \ast b)$.
        \item[(iii)]  $a \ast (b - c) =  a \ast b - a \ast c$.
        \item[(iv)]  $a \ast (b_1 + b_2 + \cdots + b_n) =  
                      a \ast b_1 + a \ast b_2 +  \cdots +  a \ast b_n$.
        \item[(v)]   \quad $a \ast (b_1 + \cdots + b_m - c_1 -\cdots - c_n) \\ =  
                     a \ast b_1 +  \cdots +  a \ast b_m - 
                     a \ast c_1 -  \cdots - a \ast c_n$.   
        \item[(vi)]  $(a \ast b + a + b) \ast c = a \ast (b \ast c) + a \ast c + b \ast c$.  
    \end{itemize}
\end{proposition}

\begin{proof}
\hfill
    \begin{itemize}
        \item [(i)] The claim follows from 0 is the identity of $(A, \circ)$ and (\ref{eqn: star}).
        \item [(ii)]  We have $0 = a \ast 0 = a \ast \big(b + (-b) \big) = a \ast b + a \ast (-b)$. Rearranging the terms gives us the claim. 
        \item[(iii)]  Since $a \ast (b - c) = a \ast b + a \ast (-c)$, the claim follows from (ii).
        \item[(iv)]  We can do mathematical induction on $n$, where inductive steps follows from (\ref{eqn: left distributive for star}).
        \item[(v)]   Rewrite the term as  $a \ast \big((b_1 + \cdots + b_m) - (c_1 + \cdots + c_n)\big)$. The claim now follows from (ii) and (iv).        
        \item[(vi)]  This follows from (\ref{eqn: star}), (iv) and basic algebraic manipulations: \\
                    \quad $(a \ast b + a + b) \ast c  \\
                                   = (a \circ b) \circ c - (a \ast b + a + b)  - c \\
                                   = a \circ (b \circ c) - a \ast b - a - b - c \\
                                   = a \ast (b \ast c + b + c) + (b \ast c + b + c)  - a \ast b - b - c \\
                                   = a \ast (b \ast c) + a \ast c + b \ast c$.\qedhere
    \end{itemize}
\end{proof}

Applying (\ref{eqn: star}) to Proposition \ref{prop: basic star arithmetic}(v), we get the following identity for the operation $\circ$. This identity will be used in the proof of Theorem \ref{main}.

\begin{proposition}
\label{prop: basic circle arithmetic}
    Let $(A, +, \circ)$ be a left brace. Then 
     \begin{align*}
         & a \circ (b_1 + \cdots + b_m - c_1 -\cdots - c_n)  \\
         = & a \circ b_1 +  \cdots +  a \circ b_m - 
                     a \circ c_1 -  \cdots - a \circ c_n +  (n-m+1)a.
     \end{align*}
   \end{proposition}

\section{Left braces with the operation $\ast$ being associative} \label{Sect: 3}
In this section, we prove Theorem~\ref{main}. 
The proof depends on Proposition \ref{prop: assoc neg star}(ii), which is obtained by
substituting (\ref{eqn: star}) into Proposition \ref{prop: assoc neg star}(i).
Proposition \ref{prop: assoc neg star}(i) is first shown in the  proof for \cite[Proposition 2.2]{CGS18}. We will sketch a proof for the convenience of the readers.

\begin{proposition}
\label{prop: assoc neg star}
    Let $(A, +, \circ)$ be a left brace such that the operation $\ast$ 
    is associative. 
    Then for all $a, b \in A$, we have
    \begin{itemize}
        \item [(i)]  $(-a) \ast b = - (a \ast b)$.
        \item [(ii)]  $(-a) \circ b = 2b - (a \circ b)$. 
    \end{itemize}
\end{proposition}

\begin{proof}
    We will sketch only the proof for (i) as (ii) can be obtained by
    substituting (\ref{eqn: star}) into (i).
    We have 
    \begin{align*}
        \big(a \ast (- a)\big) \ast b &= \big(a \ast (-a) + a + (-a)  \big) \ast b \\
                                      &= a \ast \big((-a) \ast b \big) + a \ast b + (-a) \ast b \\
                                      &= \big(a \ast (- a)\big) \ast b + a \ast b + (-a) \ast b ,
    \end{align*}
    where the second equality follows from Proposition  \ref{prop: basic star arithmetic}(vi) and the third equality follows from the
    associativity of the operation $\ast$. This implies 
    \[
        0 =  a \ast b + (-a) \ast b.
    \]
    Hence, we have $(-a) \ast b = - (a \ast b)$ for all $a, b \in A$.
\end{proof}

We now show the proof of Theorem~\ref{main}.

\begin{proof}[Proof of Theorem \ref{main}]
    Suppose the operation $\ast$ is associative. 
    Then for all $a, b, c \in A$, we have 
    \[
        (a \ast b) \ast c = a \ast (b \ast c).
    \]
    Applying (\ref{eqn: star}) twice and rearranging, we see that
    \[
    (a \circ b - a - b ) \circ c - a \circ b =  a \circ (b \circ c - b - c ) - a  -a - b \circ c + c + c.
    \]
    
    This implies 
    \[
    \begin{aligned}
        & \ \quad a^{-1} \circ \Big((a \circ b - a - b ) \circ c - a \circ b \Big) \\
        & =  a^{-1} \circ \Big(a \circ (b \circ c - b - c ) - a  -a - b \circ c + c + c  \Big).
    \end{aligned}
    \]
    Applying Proposition~\ref{prop: basic circle arithmetic},  
    substituting $a^{-1} \circ a$ with 0 and rearranging the terms gives us
    \[
          a^{-1} \circ \big(a \circ b + (- b) - a  \big) \circ c   
        = b \circ c - c  - a^{-1} \circ b \circ c +
          a^{-1} \circ c + a^{-1} \circ c.
    \]
    Applying similar manipulations, we obtain
    \[
          \big(b + a ^{-1}\circ(-b)\big) \circ c + c 
        = b \circ c + (a^{-1} \circ c + a^{-1} \circ c - a^{-1} \circ b \circ c). 
    \]
    
    Factorizing the last term on the right-hand side using Proposition~\ref{prop: basic circle arithmetic}, we have
    \[
          \big(b + a ^{-1}\circ(-b)\big) \circ c + c 
        = b \circ c + a^{-1} \circ (2c - b \circ c). 
    \]
    Proposition \ref{prop: assoc neg star}(ii) and the associativity of the operation $\circ$ gives us
    \[
    \big(b + a ^{-1}\circ(-b)\big) \circ c + c = b \circ c + a^{-1} \circ (-b) \circ c.
    \]
    Since $(A, \circ)$ is a group, for each $d \in A$, we can find the associated $a \in A$ such that $d = a^{-1} \circ (-b)$. Hence, we have for all $b,c,d \in A$,
    \[
        (b+d) \circ c + c = b \circ c + d \circ c.
    \]
    We conclude that $A$ is a two-sided brace.
\end{proof}

We remark that this proof depends on the commutativity of $(A, +)$. In fact, it is shown in \cite[\S 2.2]{KSV18} that the analog of Theorem~\ref{main} does not hold in the context of skew left braces,
where skew braces are algebraic structures which generalize left braces by relaxing the requirement that the operation + is commutative.
We refer readers to \cite{GV17} for the origin of skew left braces and 
\cite{Ven18} for a list of problems in the theory of skew left braces.

\section{Acknowledgements}

The author thanks Agata Smoktunowicz for introducing him to left braces, 
and for her encouragement in preparing this manuscript. 
The author is also grateful to Eric Jespers, Patrick Kinnear, Michael Kinyon, Leandro Vendramin and the anonymous referee for their many helpful suggestions.

\bibliographystyle{halpha}
\bibliography{main}

\end{document}